\documentclass[]{article}

\usepackage{rotating}
\usepackage{jheppub}
\usepackage{hyperref}
\usepackage{amsthm}
\usepackage{slashed}
\usepackage[utf8]{inputenc}
\usepackage[T1]{fontenc}
\usepackage{mathtools}
\usepackage[margin=1cm]{caption}

\makeatletter
\def\@seccntformat#1{%
	\expandafter\ifx\csname c@#1\endcsname\c@section\else
	\csname the#1\endcsname\quad
	\fi}
\makeatother

\def\Z{\mathbb{Z}}

\def\cF{\mathcal{F}}
\def\cG{\mathcal{G}}
\def\cH{\mathcal{H}}

\def\/{\over}

\def\t{\theta}

\def\ve{\varepsilon}

\def\d{\delta}
\def\k{\kappa}
\def\g {\gamma}

\def\G {\Gamma}

\def\Str{\mathrm{Str}}

\def\r{\mathrm}

\def\_{\hspace{2cm}}
\def\'{\:\:}

\def\be{\begin{equation}}
\def\ee{\end{equation}}
\newcommand{\bea}{\begin{eqnarray}}
\newcommand{\eea}{\end{eqnarray}}
\def\({\left(}
\def\){\right)}

\newtheorem{conjecture}{Conjecture}

\usepackage{amsmath,amsopn,amssymb,amsfonts,ctable}

\title{The Monster, the Baby Monster and Traces of Singular Moduli}

\date{2013 October 13}

\newcolumntype{R}{>{$}r<{$}}
\newcolumntype{C}{>{$}c<{$}}

\begin{document}

\begin{titlepage}
	\vfill
	\begin{flushright}
		
		
	\end{flushright}
	\vfill
	\begin{center}
		{\Large\bf The Monster, the Baby Monster and Traces of Singular Moduli}
		
		\vskip 1.4cm
		Victor Godet
		\vskip 5mm
		{\it Institute for Theoretical Physics\\
			 University of Amsterdam \\ 
			 \vskip 0.2 cm
			Science Park 904, Postbus 94485
			\\ 1090 GL Amsterdam, The Netherlands 
			}
			\vskip 0.6cm
			
				\texttt{\quad v.z.godet@uva.nl}
		\vskip 3cm
\textbf{Abstract}
\end{center}
	\vskip 0.3cm

	\begin{abstract}
		\noindent
		In this note, we provide evidence for new (super) moonshines relating the Monster and the Baby monster to some weakly holomorphic weight 1/2 modular forms defined by Zagier in his work on traces of singular moduli. They are similar in spirit to the recently discovered Thompson moonshine.
	\end{abstract}
	
	\vfill
\end{titlepage}

\section{Introduction}

Monstrous moonshine is a fascinating connection between modular forms, sporadic groups and string theory \cite{conwaynorton, norton, flm, borchmon}. Since its discovery and the explicit construction of a Monster module, many new moonshines have appeared: Mathieu moonshine \cite{eot, cheng,GaberdielI,GaberdielII,EguchiI,gannon}, generalized to Umbral moonshine \cite{dgo, um, mum}, Conway moonshine \cite{duncancon} and more recently Thompson moonshine \cite{Th, ProofTh} and O'Nan moonshine \cite{Onan}. Although their existence has been proven, most of them are still quite mysterious to this day. Furthermore, it is relevant for this work to note that Thompson and O'Nan moonshine both share a connection with traces of singular moduli \cite{zagier}.

\medskip

Here, we propose and give partial evidence for two new moonshines, which can be regarded as generalizations of Thompson moonshine for the two largest sporadic groups: the Monster and the Baby monster. Thompson moonshine involves functions defined by Zagier whose Fourier coefficients are traces of singular moduli. These functions $f_d(\tau)$ for $d\equiv 0,3\mod 4$ constitute a linear basis for the Kohnen plus-space \cite{kohnen}, which is the space of weakly holomorphic weight ${1/2}$ modular forms for $\G_0(4)$, whose Fourier coefficients $c(n)$ vanish for $n\equiv 2,3 \:\,\r{ mod }\: 4$. $f_{d}(\tau)$ is the unique function in the Kohnen plus-space whose polar part is $q^{-d}$. The modular form considered in Thompson moonshine is
\bea
\cF(\tau) &=& 2f_3(\tau)+248f_0(\tau) \\\notag
&=& 2\,q^{-3}+248+54000\, q^4 - 171990\, q^5 + 3414528\, q^8 +\dots
\eea
where $f_0(\tau) = \t(\tau) = \sum_{n\in\Z } q^{n^2}$ is  the Jacobi theta function. The Fourier coefficients of $\cF(\tau)$ have a natural decomposition in dimensions of irreducible representations of the Thompson group. Substantial evidence was given in favor of this moonshine by showing that the McKay-Thompson series $\cF_{[g]}(\tau)$, obtained by twisting $\cF(\tau)$ by an element $g  \in\mathrm{Th}$, are weight ${1/2}$ (mock) modular forms. This was done by identifying each $\cF_{[g]}(\tau)$ with a Rademacher sum projected to the Kohnen plus-space. We will follow a similar strategy to provide evidence for the moonshines proposed in this note.

\medskip

\section{Observations}

The main observations that triggered our interest are as follows. Consider the combination
\bea
\cG(\tau) &=& f_7(\tau) + f_3(\tau) - 2 f_0(\tau) \\\notag
&=& -2 + q^{-7} + q^{-3} - 4371\,q + 8315004\,q^4 - 52842475\,q^5 +5736480000 \,q^8+ \dots
\eea
We remark that the Fourier coefficients of $\cG(\tau)$ can be decomposed in dimensions of irreducible representations of the Baby monster in a simple manner.
\bea\label{ObsB}
-4371 &=& -4371 \\\notag
8315004 &=&9458750-1139374-4371-1 \\\notag\dots
\eea

A similar observation can be made for the Monster group. Let's define
\bea
\cH(\tau) &=& f_{15} (\tau)+ f_7(\tau) + f_3(\tau) - 2 f_0(\tau) \\\notag
&=& -2 + q^{-15} + q^{-7} + q^{-3} - 196884 \, q + 18517256316 \,q^4 - 292711124971 \,q^5 +\dots
\eea
and we have
\bea\label{ObsM}
-196884 &=& -196883-1 \\\notag
18517256316 &=&18538750076 - 21296876 - 196883 - 1 \\\notag\dots
\eea
where the numbers on the right are dimensions of irreducible representations of the Monster. 

\medskip

Note that the decompositions at each order are in general not unique, due to linear relations between the dimensions. In principle, the ambiguity can be resolved by requiring that all the McKay-Thompson series have good modular properties, as we will discuss later.

\medskip

This leads to the proposal of two conjectures

\begin{conjecture} There exists a graded supermodule  with a natural action of the Monster 
	\be
	W^\mathbb{M} = \bigoplus_{\substack{n\geq-15 \\ n\equiv 0,1\:\r{mod}\:4}} W^\mathbb{M} _n
	\ee
	such that the graded characters, obtained by twisting $\cH(\tau)$ by an element $g\in \mathbb{M}$ of order $|g|$, are weakly holomorphic weight ${1/2}$ (mock) modular forms for $\G_0(4|g|)$ 
	\be
	\cH_{[g]}(\tau) = \sum_{\substack{n\geq-15 \\ n\equiv 0,1\:\r{mod}\:4}} \Str_{W^\mathbb{M}_n}(g)q^n
	\ee
\end{conjecture}

\begin{conjecture} There exists a graded supermodule with a natural action of the Baby monster 
	\be
	W^\mathbb{B} = \bigoplus_{\substack{n\geq-7 \\ n\equiv 0,1\:\r{mod}\:4}} W^\mathbb{B} _n
	\ee
	such that the graded characters, obtained by twisting $\cG(\tau)$ by an element $g\in \mathbb{B}$ of order $|g|$, are weakly holomorphic weight ${1/2}$ (mock) modular forms for $\G_0(4|g|)$ 
	\be
	\cG_{[g]}(\tau) = \sum_{\substack{n\geq-7 \\ n\equiv 0,1\:\r{mod}\:4}} \Str_{W^\mathbb{B}_n}(g)q^n
	\ee
\end{conjecture}

\medskip

The two proposed moonshines appear to be very similar to Thompson moonshine.  They might be understood as generalizations obtained by adding polar terms as schematically depicted below
\be
 \underset{\makebox{\footnotesize $ $}}{1} \xrightarrow{+q^{-3}} \underset{\makebox{\footnotesize $\cF(\tau)$}}{\r{Th}}  \xrightarrow{+q^{-7}} \underset{\makebox{\footnotesize $\cG(\tau)$}}{\mathbb{B}} \xrightarrow{+q^{-15}} \underset{\makebox{\footnotesize $\cH(\tau)$}}{\mathbb{M}}
\ee
There are nonetheless notable differences worth pointing out:
\begin{itemize}
	\item $\cF(\tau)$ is related to $T_{3C}(\tau)$, the Hauptmodul associated to the Thompson group by Generalized Moonshine. The Borcherds lift
	\be
	T_{3C}(\tau) = q^{-1/3}\prod_{n>0} (1-q^n)^{A(n^2)}
	\ee
	involves the coefficients of $f_3(\tau) = q^{-3}+ \sum_{n\geq 0} A(n) q^n$. This is not sufficient to explain why the Fourier coefficients of $\cF(\tau)$ have a connection with the Thompson group since the lift only gives the $A(n)$ for $n$  a perfect square. 
	
	For the moonshines proposed here, no such link appears to exist with Generalized Moonshine. $\cH(\tau)$ seems unrelated to $T_{1A}(\tau) = j(\tau)$ and $\cG(\tau)$ seems unrelated to $T_{2A}(\tau)$, where $T_{1A}(\tau)$ and $T_{2A}(\tau)$ are the Hauptmoduls associated with the Monster and the Baby monster by Generalized Moonshine.
	
	\item The Thomspon supermodule has a simple graded structure
	\be
	W^{\r{Th}} = \sum_{\substack{n\geq 3 \\ n\equiv 0,1\: \r{ mod } \:4}}^\infty W^{\r{Th}}_n
	\ee
	where the $\Z_2$-grading is such that $W^{\r{Th}}_n$ is purely even for $n\equiv 0 \,\:\mathrm{mod }\: 4$ and purely odd for $n\equiv 1 \,\:\r{mod }\: 4$, except for $W^{\r{Th}}_{-3}$ which is purely even. 
	
	In our cases, the supermodules will have generically both an even and an odd part at each $n$.
\end{itemize}

\section{Evidence}

We will now give evidence in favor of these two conjectures. This is done by investigating the modular properties of the McKay-Thompson series $\cH_{[g]}(\tau)$ and $\cG_{[g]}(\tau) $ for each conjugacy classes $[g]$ of the Monster and the Baby monster respectively. We will show that many of these $q$-series can be written naturally as certain Rademacher sums \cite{ChengDuncanI,ChengDuncanII,DuncanFrenkel}. Hence, they seem to be weakly holomorphic weight ${1/2}$ (mock) modular forms for $\G_0(4|g|)$, where $|g|$ is the order of the conjugacy class. If the observations $(\ref{ObsB})$ and $(\ref{ObsM})$ were coincidental or trivial, we would not expect the McKay-Thompson series to have any modular property at all. This gives some strong support for the moonshines proposed in this note.

\medskip

In Thompson moonshine, the McKay-Thompson series have been shown to be of the following form
\be
\cF_{[g]}(\tau) = 2 Z_{N_g,\psi_g}(\tau)+ \sum_{\substack{m>0\\ m^2 | (h_g|g|)} }{\k_{m,g} \t(m^2\tau)}
\ee
where $Z_{N_g,\psi_g}(\tau)$ is a Rademacher sum  projected to the Kohnen plus-space. Each of these functions is fixed by the following data
\begin{itemize}
	\item $N_g$ specifies the associated modular group $\G_0(4N_g)$. Here, $N_g =|g|$ where $|g|$ is the order of $g\in \r{Th}$.
	\item $\psi_g$ is a multiplier system which is specified by two integers $(v_g,h_g)$ and
	\be
	\psi_g(\g) = \psi_0(\g) e\left(-v_g {cd\/N_g h_g} \right) \qquad  \g = \begin{pmatrix} a & b \\ c & d\end{pmatrix} \in \G_0(4|g|)
	\ee
	where $e(z) = e^{2\pi i z}$ and $\psi_0(\g) = \left(c\/d \right) \ve_d$ is the multiplier system associated with the Jacobi theta function $\t(\tau) =\sum_{n\in\Z}q^{n^2}$ 
	\item $\k_{m,g}$ are rational numbers which specify theta-corrections
\end{itemize}

\medskip

For the Monster, we analogously want to identify the McKay-Thompson series $\cH_{[g]}(\tau)$, where $g\in\mathbb{M}$, with the following ansatz
\be\label{RadF}
 Z^{[-15]}_{N^{(15)}_g,\psi_g^{(15)}}(\tau)+ Z^{[-7]}_{N^{(7)}_g,\psi_g^{(7)}}(\tau)+ Z^{[-3]}_{N^{(3)}_g,\psi_g^{(3)}}(\tau) + \sum_{\substack{m>0\\ m^2 | (h_g |g|)} }{\k_{m,g} \t(m^2\tau)}
\ee
This ansatz is natural from the polar structure of $\cH(\tau) = q^{-15}+q^{-7} +q^{-3}+\dots$. The following data need to be specified
\begin{itemize}
	\item $N^{(15)}_g, N^{(7)}_g,N^{(3)}_g$ are integers such that $Z^{[-\mu]}_{N_g^{(\mu)}}$, where $\mu=15, 7, 3$, is a (mock) modular form for $\G_0(4N_g^{(\mu)})$ and
	\be
	\mathrm{lcm} (N_g^{(15)},N_g^{(7)} ,N_g^{(3)}) = |g|
	\ee
	\item   $\psi_g^{(15)},\psi_g^{(7)},\psi_g^{(3)}$ are the multiplier systems specified by two integers $(v_g^{(\mu)},h_g^{(\mu)})$ and
	\be
	\psi_g^{(\mu)}(\g) = \psi_0(\g) e\left(-v_g^{(\mu)} {cd\/N_g^{(\mu)} h_g^{(\mu)}} \right)\qquad  \g = \begin{pmatrix} a & b \\ c & d\end{pmatrix}\in \G_0(4|g|)
	\ee
	and we also define $h_g =  \mathrm{lcm} (h_g^{(15)},h_g^{(7)} ,h_g^{(3)})$
	\item $\k_{m,g}$ are rational numbers which specify theta-corrections 

\end{itemize}
Given this data, the function defined in (\ref{RadF}) is a (mock) modular form of weight $1/2$ for $\G_0(4|g|)$. 

\medskip

For the Baby monster, we want to identify the McKay-Thompson series $\cG_{[g]}(\tau)$, for $g\in\mathbb{B}$, with the following ansatz

\be\label{RadFB}
  Z^{[-7]}_{N^{(7)}_g,\psi_g^{(7)}}(\tau)+ Z^{[-3]}_{N^{(3)}_g,\psi_g^{(3)}}(\tau) + \sum_{\substack{m>0\\ m^2 | (h_g|g|)} }{\k_{m,g} \t(m^2\tau)}
\ee
where we use analogous definitions.

\medskip

The Rademacher sums projected to Kohnen plus-space have Fourier expansion
\be
Z^{[-m]}_{N,\psi}(\tau) = q^{-m} + \sum_{\substack{n\geq 0\\ n\equiv 0,1\,\r{ mod }\,4}} c_N(-m,n)\,q^n
\ee
and the explicit formulas for their Fourier coefficients are given in \cite{mp} 
\bea
n=0&:\quad&c_N(-m;0) = 4\,\pi \sqrt{m}(1-i) \sum_{\substack{c>0\\4N|c}} (1+\d_{\text{odd}}(c/4)) {K_\psi(-m,0,c)\/c^{3/2}}\\\notag
n>0 &:\quad&c_N(-m;n) = \pi\sqrt{2} \left(n\/m \right)^{-1/4} (1-i) \sum_{\substack{c>0\\4N|c}} (1+\d_{\text{odd}}(c/4)) {K_\psi(-m,n,c)\/c} I_{1\/2} \left(4\pi\sqrt{mn}\/c \right)
\eea
where
\be
K_\psi(m,n,c) = \sum_{d \in (\Z/c\Z)^*} \psi(c,d)\left(c\/d \right) \ve_d e\left(m \bar{d}+nd\/c \right),  \qquad \bar{d} d \equiv 1 \mod c
\ee
and
\be
\ve_d = \begin{cases} 1  &d\equiv 1\mod 4\\ i   &d\equiv 3\mod 4 \end{cases}  \qquad \qquad  \d_{\text{odd}}(k) = \begin{cases} 1 & k \text{ odd} \\ 0   & k \text{ even}\end{cases}
\ee
With these definitions, $ Z^{[-m]}_{N,\psi}(\tau)  $ is a weakly holomorphic mock modular form of weight ${1/2}$ for $\G_0(4N)$. These formulas have been derived for $N$ odd but it can be argued that they are also valid for $N$ even, when $\G_0(4N)$ is genus zero \cite{Th}. The even $N$ cases that we have checked fall into this category \cite{genuszero}.

\medskip

Now, we will try to identify the McKay-Thompson series with Rademacher sums. Of course, different decompositions are going to give different McKay-Thompson series.  Since we are dealing with a supermodule, we can add any linear relation between the dimensions to a given decomposition.  For example, we have the following relation between the dimensions of irreducible representations of the Baby monster
\be
9550635 - 9458750 - 96255 + 4371 - 1 = 0
\ee
Adding this in the decompositions will only change a few number of classes and it is rather difficult to check all the possibilities, especially for terms of high $q$ power. Among all the possible "reasonable" decompositions, we chose the one that matches with our ansatz for the highest number of classes, among those that we have been able to check.

\medskip

Using this strategy, we have settled to the following decompositions for the first few coefficients. For the Monster,
\bea
-196884 &=& -196883-1 \\\notag
18517256316 &=&19360062527 - 842609326 - 196883 - 2 \\\notag
-292711124971 &=&-293553734298+ 842609326+1 \\\notag
312217411718400&=&190292345709543 + 125510727015275- 3879214937598+ 293553734298 \\\notag && + 196883 - 1 \\\notag
-2374124840259859 &=& -2374124840062976-196883 
\eea
and for the Baby monster,
\bea
-4371 &=& -4371 \\\notag
8315004 &=& 9550635-1139374-96255-2 \\\notag
-52842475 &=&-63532485 + 9458750 + 1139374+ 96255 - 4371 + 2 \\\notag
5736480000 &=&  4275362520 + 3214743741 - 1407126890 - 356054375 + 9550635 + 4371 - 2
\eea
This translates in the following decompositions of the supermodule $W^{\mathbb{M}}$ into irreducible representations of the Monster. We denote them by $V_1,V_2,\dots$ with $V_1$ being the trivial representation.
\bea\label{decM}
W^{\mathbb{M}}_1 &=& V_2^-\oplus V_1^- \\\notag
W^{\mathbb{M}}_2 &=& V_6^+ \oplus V_4^- \oplus V_2^-\oplus \,2\cdot V_1^- \\\notag
W^{\mathbb{M}}_3 &=& V_4^+ \oplus V_1^+ \oplus V_7^- \\\notag
W^{\mathbb{M}}_4 &=& V_{11}^+\oplus V_{10}^+ \oplus V_7^+ \oplus V_2^+ \oplus V_8^- \oplus V_1^- \\\notag
W^{\mathbb{M}}_5 &=& V_{15}^- \oplus V_{2}^-
\eea
and the superscript $\pm$ refers to the $\Z_2$-grading of the supermodule.

\medskip

Similarly, we denote by $U_1,U_2,\dots$ the irreducible representations of the Baby monster starting from the trivial one and we have
\bea\label{decB}
W^{\mathbb{B}}_1 &=& U_2^-\\\notag
W^{\mathbb{B}}_2 &=& U_6^+ \oplus U_4^- \oplus U_3^-\oplus \,2\cdot U_1^- \\\notag
W^{\mathbb{B}}_3 &=& U_3^+\oplus U_5^+\oplus U_4^+ \oplus \,2\cdot U_1^+ \oplus U_7^-\oplus U_2^- \\\notag
W^{\mathbb{B}}_4 &=& U_{13}^+ \oplus U_{11}^+ \oplus U_6^+ \oplus U_2^+ \oplus U_{10}^- \oplus U_9^- \oplus \,2\cdot U_1^-
\eea
We have computed all the McKay-Thompson series based on these decompositions. They are tabulated in the Appendix.
It is rather hard to go to higher order in $q$ since the number of possible decompositions grows quickly.

\medskip
Furthermore, we have restricted the possible multiplier systems by defining
\be
\hat{h} = {h\/ \r{gcd}(h,4)} \qquad \hat{v} = {4v \/ \r{gcd}(h,4)} \mod \hat{h}
\ee
and imposing $\hat{v} = \pm 1 \mod \hat{h}$.  This contains all the cases that appear in Thompson moonshine \cite{ProofTh}. In practice, we don't seem to get new identifications by relaxing this condition.

\medskip

We will now give some examples of the identifications that we have found. A more extensive list can be found in the Appendix. For the Monster, some McKay-Thompson series are
\bea
\cH_{1A}(\tau)&=& -2 + q^{-15} + q^{-7} + q^{-3}-196884\,q+18517256316\,q^4 \\\notag
&&\quad -292711124971\,q^5+312217411718400\,q^8-2374124840259859\,q^9+\dots\\
\cH_{2A}(\tau) &=& -2 + q^{-15} + q^{-7} + q^{-3} -4372\,q+8218748\,q^4\\\notag&& \quad-52842475\,q^5+5746986240\,q^8 -22509166867\,q^9+\dots \\
\cH_{2B}(\tau) &=&-2 + q^{-15} + q^{-7} + q^{-3} -276\,q-71556\,q^4-85995\,q^5\\\notag&&\quad+12258560\,q^8-4100371\,q^9+\dots \\
\cH_{3B}(\tau) &=&-2 + q^{-15} + q^{-7} + q^{-3}-54\,q+1674\,q^4+3707\,q^5-33150\,q^8 +1403\,q^9+\dots 
\eea
and can be nicely identified with our ansatz (\ref{RadF})
\bea
\cH_{1A}(\tau)&=&
Z^{[-15]}_{1}(\tau)+ Z^{[-7]}_{1}(\tau)+ Z^{[-3]}_{1}(\tau) -42\,\t(\tau)  \\
\cH_{2A}(\tau) &=& Z^{[-15]}_{2}(\tau)+ Z^{[-7]}_{1}(\tau)+ Z^{[-3]}_{1}(\tau) -26\,\t(\tau)  \\
\cH_{2B}(\tau) &=& Z^{[-15]}_{2}(\tau)+ Z^{[-7]}_{2}(\tau)+ Z^{[-3]}_{1}(\tau) -18\,\t(\tau) \\
\cH_{3B}(\tau) &=& Z^{[-15]}_{3}(\tau)+ Z^{[-7]}_{3}(\tau)+ Z^{[-3]}_{3}(\tau) +{15\/2}\,\t(\tau)
\eea
which we have checked up to $q^9$. From now on, we will skip the multiplier subscript $\psi$ in the Rademacher sums $Z^{[-m]}_N(\tau)$, whenever it is equal to $\psi_0$.

\medskip

For the Baby monster, some McKay-Thompson series are
\bea
\cG_{1A}(\tau)&=&-2+q^{-7}+q^{-3} -4371\,q+8315004\,q^4-52842475\,q^5+5736480000\,q^8+\dots  \\
\cG_{2B}(\tau) &=&  -2+q^{-7}+q^{-3} -275\,q+24700\,q^4-85995\,q^5+1752320\,q^8+\dots \\
\cG_{2D}(\tau) &=&  -2+q^{-7}+q^{-3} -19\,q-1924\,q^4+21\,q^5+44288\,q^8+\dots  \\
\cG_{3B}(\tau) &=&  -2+q^{-7}+q^{-3} +3\,q+30\,q^4-181\,q^5+870\,q^8+\dots   \\
\cG_{4G}(\tau) &=&  -2+q^{-7}+q^{-3} -11\,q-132\,q^4+21\,q^5+768\,q^8+\dots
\eea
and are identified as
\bea
\cG_{1A}(\tau)&=&Z^{[-7]}_{1}(\tau)+ Z^{[-3]}_{1}(\tau) -18\,\t(\tau)  \\
\cG_{2B}(\tau) &=&  Z^{[-7]}_{2}(\tau)+ Z^{[-3]}_{1}(\tau) -10\,\t(\tau)  \\
\cG_{2D}(\tau) &=&  Z^{[-7]}_{2}(\tau)+ Z^{[-3]}_{2}(\tau) -2\,\t(\tau)  \\
\cG_{3B}(\tau) &=&  Z^{[-7]}_{3}(\tau)+ Z^{[-3]}_{3}(\tau) +{9\/2}\,\t(\tau)   \\
\cG_{4G}(\tau) &=&  Z^{[-7]}_{4}(\tau)+ Z^{[-3]}_{2,(1,8)}(\tau) +2\,\t(\tau)  -8\,\t(4\tau)
\eea
which we have checked up to $q^8$. Here, $Z^{[-3]}_{2,(1,8)}$ is a notation for $Z^{[-3]}_{2,\psi}$ where $\psi$ is the multiplier system defined by $(v,h)=(1,8)$. More such identifications are listed in the Appendix.

\medskip

In this note, we have tested all the conjugacy classes $[g]$ up to order $|g|= 11$. For the Monster, the McKay-Thompson series can be identified as Rademacher sums of the form (\ref{RadF}) for the classes
\be
1A, 2A, 2B, 3B, 4A, 5B, 6D, 6E, 7B, 8B,  8D, 8E, 10C,  10E
\ee 
but not for the following ones 
\be
3A, 3C, 4B, 4D, 5A, 6A, 6B, 6C, 6F, 7A, 8A,	8C, 8F, 9A, 9B, 10A, 10B, 10D, 11A 
\ee

For the Baby monster,  the McKay-Thompson series can be identified as Rademacher sums of the form (\ref{RadFB}) for the classes
\be
1A, 2B, 2D, 4G,  5B, 6G, 6J, 7A, 8D, 8F, 8H, 8I, 9A, 9B, 10D, 10F
\ee
but not for the following ones 
\bea
2A, 2C, 3A, 4A, 4B, 4C, 4D, 4E, 4F, 4H, 4I, 4J, 5A, 6A, 6B, 6C, 6D, 6E, 6F, \\
6H, 6I, 6K, 8A, 8B, 8C, 8E,  8G, 8J, 8K, 8L, 8M, 8N, 10A, 10B, 10C, 10E, 11A
\eea

The fact that many conjugacy classes cannot be matched suggests that the ansatz (\ref{RadF}), as well as (\ref{RadFB}), is too restrictive. For example, for the class $2A$ of the Baby monster 
\be
\cG_{2A}(\tau) = -2+q^{-7}+q^{-3} + 493\,q + 146300\,q^4 +565269\,q^5 +  18517248\,q^8 + \dots
\ee
which cannot be identified with a Rademacher sum of the form (\ref{RadFB}). However, it is nicely matched with the following Rademacher sum
\be
\cG_{2A}(\tau) = Z^{[-7]}_{2}(\tau)+ Z^{[-3]}_{2}(\tau) +Z^{[-4]}_{1}(\tau) - Z^{[-4]}_{2}(\tau) - 10\,\t(\tau)
\ee
This suggests that we should consider more general combinations of Rademacher sums involving other polar terms. In the above example, we have some Rademacher sums with polar term $q^{-4}$ and the combination that appears is
\be \label{4term}
Z^{[-4]}_{1}(\tau) - Z^{[-4]}_{2}(\tau)
\ee
Note that the difference is needed in order to cancel the unwanted $q^{-4}$.

\medskip

More generally, it seems that we need to consider linear combinations of the Rademacher sums $Z^{[-m]}(\tau)$ for $-15 \leq m \leq -3$, with different levels and multiplier systems, such that the polar part of the combination is $q^{-15}+q^{-7}+q^{-3}$ for the Monster. Similarly,  we will have $-7\leq m\leq -3$ for the Baby monster such that the polar part is $q^{-7}+q^{-3}$. However, this allows for a very large set of functions, which is more difficult to handle. Since we have demonstrated several identifications, we expect the general ansatz to accommodate all conjugacy classes.

\section{Conclusion}

We have proposed two new moonshines relating the Monster and the Baby monster to some weakly holomorphic weight $1/2$ modular forms. We have provided some evidence, demonstrating that many McKay-Thompson series can be identified as Rademacher sums. Despite the similarity to Thompson moonshine, we have also seen that a more general ansatz for the McKay-Thompson series is needed in order to accommodate all classes of the Monster and the Baby monster. Even though the evidence we provide is not complete, it suggests that all conjugacy classes should be identifiable with this more general ansatz. 

\medskip

We also note that the proposed moonshine for the Monster seems to be unrelated to the familiar Monstrous moonshine. It involves a supermodule with characters of weight ${1/2}$ rather than a module with characters of weight $0$. 

\medskip

It would be highly interesting to find a natural construction of these supermodules, perhaps from a super vertex operator algebra. 

\vskip 0.5cm

\paragraph{Acknowledgments}

\medskip

I would like to thank Francesca Ferrari for useful discussions, and especially Vassilis Anagiannis for his valuable help in most stages of this project.

\newpage
\appendix

\section{Identifications with Rademacher sums}

This is a complete list for $|g|\leq 11$ of the proposed identifications between the McKay-Thompson series $\cG_{[g]}(\tau)$ of the Baby monster with Rademacher sums of the form (\ref{RadFB}). Note that the McKay-Thompson series are computed assuming the decomposition (\ref{decB}). $Z^{[-m]}_{N,(v,h)}$ is a notation for $Z^{[-m]}_{N,\psi}$ where $\psi$ is the multiplier system defined by the integers $(v,h)$.
\bea
\cG_{1A}(\tau)&=&Z^{[-7]}_{1}(\tau)+ Z^{[-3]}_{1}(\tau) -18\,\t(\tau)  \\
\cG_{2B}(\tau) &=&  Z^{[-7]}_{2}(\tau)+ Z^{[-3]}_{1}(\tau) -10\,\t(\tau)  \\
\cG_{2D}(\tau) &=&  Z^{[-7]}_{2}(\tau)+ Z^{[-3]}_{2}(\tau) -2\,\t(\tau)  \\
\cG_{3B}(\tau) &=&  Z^{[-7]}_{3}(\tau)+ Z^{[-3]}_{3}(\tau) +{9\/2}\,\t(\tau)   \\
\cG_{4G}(\tau) &=&  Z^{[-7]}_{4}(\tau)+ Z^{[-3]}_{2,(1,8)}(\tau) +2\,\t(\tau)  -8\,\t(4\tau)\\
\cG_{5B}(\tau) &=&  Z^{[-7]}_{5}(\tau)+ Z^{[-3]}_{5}(\tau) +2\,\t(\tau)   \\
\cG_{6G}(\tau) &=&  Z^{[-7]}_{6}(\tau)+ Z^{[-3]}_{3}(\tau) +{1\/2}\,\t(\tau)\\
\cG_{6J}(\tau) &=&  Z^{[-7]}_{6}(\tau)+ Z^{[-3]}_{6}(\tau) -{1\/2}\,\t(\tau)\\
\cG_{7A}(\tau) &=&  Z^{[-7]}_{7}(\tau)+ Z^{[-3]}_{7}(\tau)  -{5\/2}\,\t(\tau)   \\
\cG_{8D}(\tau) &=&  Z^{[-7]}_{8}(\tau)+ Z^{[-3]}_{4}(\tau)  -4\,\t(\tau)   +2\,\t(4\tau)\\
\cG_{8F}(\tau) &=&  Z^{[-7]}_{8}(\tau)+ Z^{[-3]}_{4}(\tau)  +4\,\t(\tau)   -6\,\t(4\tau)\\
\cG_{8H}(\tau) &=&  Z^{[-7]}_{8}(\tau)+ Z^{[-3]}_{4}(\tau)  -2\,\t(4\tau)\\
\cG_{8I}(\tau) &=&  Z^{[-7]}_{8,(1,16)}(\tau)+ Z^{[-3]}_{4,(1,8)}(\tau)  +2\,\t(\tau) + 4\,\t(4\tau)-8\,\t(16\tau)\\
\cG_{9A}(\tau) &=&  Z^{[-7]}_{9}(\tau)+ Z^{[-3]}_{9}(\tau)  +{3\/2}\,\t(\tau) -3\,\t(9\tau)   \\
\cG_{9B}(\tau) &=&  Z^{[-7]}_{9}(\tau)+ Z^{[-3]}_{9}(\tau)  -3\,\t(\tau) +{3\/2}\,\t(9\tau)   \\
\cG_{10D}(\tau) &=&  Z^{[-7]}_{10}(\tau)+ Z^{[-3]}_{5}(\tau) \\
\cG_{10F}(\tau) &=&  Z^{[-7]}_{10}(\tau)+ Z^{[-3]}_{10}(\tau)  -2\,\t(\tau)
\eea

\newpage
These are the McKay-Thompson series $\cH_{[g]}$ of the Monster, with $|g|\leq 11$, that can be matched with Rademacher sums of the form (\ref{RadF}). The McKay-Thompson series are computed assuming the decomposition (\ref{decM}).
\bea
\cH_{1A}(\tau)&=&Z^{[-15]}_{1}(\tau)+ Z^{[-7]}_{1}(\tau)+ Z^{[-3]}_{1}(\tau) -42\,\t(\tau)  \\
\cH_{2A}(\tau) &=& Z^{[-15]}_{2}(\tau)+ Z^{[-7]}_{1}(\tau)+ Z^{[-3]}_{1}(\tau) -26\,\t(\tau)  \\
\cH_{2B}(\tau) &=& Z^{[-15]}_{2}(\tau)+ Z^{[-7]}_{2}(\tau)+ Z^{[-3]}_{1}(\tau) -18\,\t(\tau) \\
\cH_{3B}(\tau) &=& Z^{[-15]}_{3}(\tau)+ Z^{[-7]}_{3}(\tau)+ Z^{[-3]}_{3}(\tau) +{15\/2}\,\t(\tau)\\
\cH_{4A}(\tau) &=& Z^{[-15]}_{4}(\tau)+ Z^{[-7]}_{4}(\tau)+ Z^{[-3]}_{1}(\tau) -18\,\t(\tau)+12\,\t(4\tau) \\
\cH_{4C}(\tau) &=& Z^{[-15]}_{4}(\tau)+ Z^{[-7]}_{4}(\tau)+ Z^{[-3]}_{2}(\tau) -10\,\t(\tau)+12\,\t(4\tau) \\
\cH_{5B}(\tau) &=& Z^{[-15]}_{5}(\tau)+ Z^{[-7]}_{5}(\tau)+ Z^{[-3]}_{5}(\tau) +3\,\t(\tau)\\
\cH_{6D}(\tau) &=& Z^{[-15]}_{6}(\tau)+ Z^{[-7]}_{3}(\tau)+ Z^{[-3]}_{3}(\tau) +{11\/2}\,\t(\tau)\\
\cH_{6E}(\tau) &=& Z^{[-15]}_{6}(\tau)+ Z^{[-7]}_{6}(\tau)+ Z^{[-3]}_{3}(\tau) +{3\/2}\,\t(\tau)\\
\cH_{7B}(\tau) &=& Z^{[-15]}_{7}(\tau)+ Z^{[-7]}_{7}(\tau)+ Z^{[-3]}_{7}(\tau) +{3\/2}\,\t(\tau)\\
\cH_{8B}(\tau) &=& Z^{[-15]}_{8}(\tau)+ Z^{[-7]}_{8}(\tau)+ Z^{[-3]}_{2,(1,8)}(\tau) -6\,\t(4\tau)\\
\cH_{8D}(\tau) &=& Z^{[-15]}_{8}(\tau)+ Z^{[-7]}_{8}(\tau)+ Z^{[-3]}_{4}(\tau) +8\,\t(\tau)-10\,\t(4\tau)\\
\cH_{8E}(\tau) &=& Z^{[-15]}_{8}(\tau)+ Z^{[-7]}_{8}(\tau)+ Z^{[-3]}_{4}(\tau) +4\,\t(\tau)-6\,\t(4\tau)\\
\cH_{10C}(\tau) &=& Z^{[-15]}_{10}(\tau)+ Z^{[-7]}_{5}(\tau)+ Z^{[-3]}_{5}(\tau) +{5\/2}\,\t(\tau)\\
\cH_{10E}(\tau) &=& Z^{[-15]}_{10}(\tau)+ Z^{[-7]}_{10}(\tau)+ Z^{[-3]}_{5}(\tau) +{1\/2}\,\t(\tau)
\eea

\section{Tables of McKay-Thompson series}

We provide here the Fourier coefficients of all the McKay--Thompson series $\cG_{[g]}(\tau)$ up to order $q^{8}$, and $\cH_{[g]}(\tau)$ up to order $q^{9}$. The character tables for the Monster and the Baby monster were obtained with the GAP4 computer algebra package \cite{GAP4}.

\begin{sidewaystable}\setlength{\tabcolsep}{-1pt}\renewcommand{\arraystretch}{1}
	\begin{small}
		\centering
		\caption{Coefficients of $q^n$ in the McKay--Thompson series $\cH_{[g]}(\tau)$ for the Monster, part one.}\label{mtM}\smallskip
		\resizebox{!}{64pt}{\begin{tabular}
			{C@{ }@{\;}R@{ }R@{ }R@{ }R@{ }R@{ }R@{ }R@{ }R@{ }R@{ }R@{ }R@{ }R@{ }R@{ }R@{ }R@{ }R@{ }R@{ }R@{ }R@{ }R@{ }R@{ }R@{ }R@{ }R@{ }R@{ }R}\toprule
			 n ~\backslash ~[g]&
			1A & 2A & 2B & 3A & 3B & 3C & 4A & 4B & 4C & 4D & 5A & 5B & 6A & 6B & 6C & 6D & 6E & 6F  \\
			\midrule
			-15 &1 & 1 & 1 & 1 & 1 & 1 & 1 & 1 & 1 & 1 & 1 & 1 & 1 & 1 & 1 & 1 & 1 & 1 \\
			-7 & 1 & 1 & 1 & 1 & 1 & 1 & 1 & 1 & 1 & 1 & 1 & 1 & 1 & 1 & 1 & 1 & 1 & 1  \\
			-3	& 1 & 1 & 1 & 1 & 1 & 1 & 1 & 1 & 1 & 1 & 1 & 1 & 1 & 1 & 1 & 1 & 1 & 1 	\\
			0     &-2 & -2 & -2 & -2 & -2 & -2 & -2 & -2 & -2 & -2 & -2 & -2 & -2 & -2 & -2 & -2 & -2 & -2  \\
			1 & -196884 & -4372 & -276 & -783 & -54 & 0 & -276 & -52 & -20 & 12 & -134 & -9 & -79 & -78 & -15 & 2 & -6 & 0  \\
			4  &18517256316 & 8218748 & -71556 & 240786 & 1674 & 0 & 26748 & -1668 & 124 & -132 & 5566 & 66 & 1298 & 1650 & -366 & 2 & -6 & 0\\
			5  & -292711124971 & -52842475 & -85995 & -999397 & 3707 & -4123 & -85995 & -491 & 21 & 533 & -14971 & 29 & -4069 & -3861 & 27 & -181 & 27 & -27  \\
			8  & 312217411718400 & 5746986240 & 12258560 & 38137290 & -33150 & 34998 & 1707264 & 46848 & -768 & -3328 & 187150 & -350 & 38346 & 33770 & 3530 & 762 & -142 & 182  \\
			9 & -2374124840259859 & -22509166867 & -4100371 & -110509894 & 1403 & 61505 & -4100371 & -8243 & -4115 & 13 & -393109 & -484 & -62854 & -64141 & -70 & 1379 & 11 & 65 \\
			\end{tabular}}
			\resizebox{!}{66pt}{
				\begin{tabular}{C@{ }@{\;}R@{ }R@{ }R@{ }R@{ }R@{ }R@{ }R@{ }R@{ }R@{ }R@{ }R@{ }R@{ }R@{ }R@{ }R@{ }R@{ }R@{ }R@{ }R@{ }R@{ }R@{ }R@{ }R@{ }R@{ }R@{ }R@{ }R@{ }R@{ }R@{ }R@{ }R@{ }R@{ }R@{ }R@{ }R@{ }R@{ }R@{ }R}\midrule
					n ~\backslash ~[g] &7A & 7B & 8A & 8B & 8C & 8D & 8E & 8F & 9A & 9B & 10A & 10B & 10C & 10D & 10E & 11A & 12A & 12B & 12C & 12D & 12E & 12F & 12G & 12H & 12I & 12J & 13A & 13B & 14A & 14B & 14C\\
			\midrule
			-15 &1 & 1 & 1 & 1 & 1 & 1 & 1 & 1 & 1 & 1 & 1 & 1 & 1 & 1 & 1 & 1 & 1 & 1 & 1 & 1 & 1 & 1 & 1 & 1 & 1 & 1 & 1 & 1 & 1 & 1& 1  \\
			-7 & 1 & 1 & 1 & 1 & 1 & 1 & 1 & 1 & 1 & 1 & 1 & 1 & 1 & 1 & 1 & 1 & 1 & 1 & 1 & 1 & 1 & 1 & 1 & 1 & 1 & 1 & 1 & 1 & 1 & 1& 1   \\
			-3	& 1 & 1 & 1 & 1 & 1 & 1 & 1 & 1 & 1 & 1 & 1 & 1 & 1 & 1 & 1 & 1 & 1 & 1 & 1 & 1 & 1 & 1 & 1 & 1 & 1 & 1 & 1 & 1 & 1 & 1 & 1 	\\
			0     &-2 & -2 & -2 & -2 & -2 & -2 & -2 & -2 & -2 & -2 & -2 & -2 & -2 & -2 & -2 & -2 & -2 & -2 & -2 & -2 & -2 & -2 & -2 & -2 & -2 & -2 & -2 & -2 & -2 & -2& -2    \\
			1 & -51 & -2 & -36 & -12 & 0 & 4 & -4 & 0 & -27 & 0 & -22 & -6 & 3 & -21 & -1 & -17 & -15 & -6 & -7 & 0 & 1 & -6 & 2 & -14 & -2 & 0 & -12 & 1 & -11 & -3 & -10 \\
			4  &594 & 6 & 252 & -132 & -52 & -4 & -4 & 12 & 108 & 0 & -2 & -66 & -2 & 54 & -6 & 20 & 18 & 18 & -30 & 0 & -14 & -42 & 6 & 10 & -2 & 0 & 0 & 0 & -22 & -30 & -2  \\
			5  & -1273 & 1 & -491 & 21 & 53 & 21 & 21 & -11 & -208 & 8 & -75 & 5 & 25 & -95 & 5 & -45 & 27 & 27 & -5 & -27 & 27 & 11 & -5 & -21 & 3 & 5 & -5 & 8 & 7 & 7 & -7 \\
			8  & 8448 & 20 & 2560 & 768 & -160 & 0 & 0 & -32 & 924 & 6 & 590 & 270 & -10 & 350 & 10 & 166 & -54 & -54 & 138 & 54 & -54 & 110 & 30 & 66 & 6 & -10 & 32 & 6 & 72 & 48 & 20  \\
			9 & -14603 & -1 & -4131 & -11 & 1 & 5 & -3 & 1 & -1351 & 26 & -517 & -21 & 8 & -496 & 4 & -192 & -70 & 11 & -62 & 65 & 10 & -5 & 19 & -77 & 7 & 1 & -37 & 2 & -19 & 5 & -9 \\
			\bottomrule
		\end{tabular}}
	\end{small}
\end{sidewaystable}

\begin{sidewaystable}
	\begin{small}
		\centering
		\caption{Coefficients of $q^n$ in the McKay--Thompson series $\cH_{[g]}(\tau)$ for the Monster, part two.}\smallskip
		\begin{tabular}{C@{ }@{\;}R@{ }R@{ }R@{ }R@{ }R@{ }R@{ }R@{ }R@{ }R@{ }R@{ }R@{ }R@{ }R@{ }R@{ }R@{ }R@{ }R@{ }R@{ }R@{ }R@{ }R@{ }R@{ }R@{ }R@{ }R@{ }R@{ }R@{ }R@{ }R@{ }R@{ }R@{ }R@{ }R@{ }R@{ }R@{ }R@{ }R}\toprule
			n ~\backslash ~[g] &15A & 15B & 15C & 15D & 16A & 16B & 16C & 17A & 18A & 18B & 18C & 18D & 18E & 19A & 20A & 20B & 20C & 20D & 20E & 20F & 21A & 21B & 21C & 21D & 22A & 22B & 23A & 23B & 24A \\
			\midrule
			-15 &1 & 1 & 1 & 1 & 1 & 1 & 1 & 1 & 1 & 1 & 1 & 1 & 1 & 1 & 1 & 1 & 1 & 1 & 1 & 1 & 1 & 1 & 1 & 1 & 1 & 1 & 1 & 1& 1   \\
			-7 & 1 & 1 & 1 & 1 & 1 & 1 & 1 & 1 & 1 & 1 & 1 & 1 & 1 & 1 & 1 & 1 & 1 & 1 & 1 & 1 & 1 & 1 & 1 & 1 & 1 & 1 & 1 & 1 & 1  \\
			-3	& 1 & 1 & 1 & 1 & 1 & 1 & 1 & 1 & 1 & 1 & 1 & 1 & 1 & 1 & 1 & 1 & 1 & 1 & 1 & 1 & 1 & 1 & 1 & 1 & 1 & 1 & 1 & 1  & 1 	\\
			0     &-2 & -2 & -2 & -2 & -2 & -2 & -2 & -2 & -2 & -2 & -2 & -2 & -2 & -2 & -2 & -2 & -2 & -2 & -2 & -2 & -2 & -2 & -2 & -2 & -2 & -2 & -2 & -2 & -2     \\
			1 & -8 & 1 & -9 & 0 & -4 & 0 & -8 & -7 & 2 & -7 & -3 & 0 & -6 & -6 & -6 & -2 & -1 & 2 & -3 & -5 & -6 & 1 & 0 & -5 & -5 & -1 & -4 & -4 & -3  \\
			4  & -14 & 4 & -6 & 0 & -20 & -4 & -4 & -8 & 2 & -16 & -12 & 0 & -6 & -6 & -2 & -18 & -2 & -2 & -12 & -6 & 0 & 0 & 0 & -6 & -12 & -12 & -7 & -7 & -6 \\
			5  & -22 & -13 & 2 & 2 & 5 & 5 & 5 & 2 & 8 & 8 & 0 & 0 & 0 & 6 & 5 & 9 & 5 & 13 & 3 & 1 & 14 & 7 & -7 & 4 & 7 & 3 & 2 & 2 & 3 \\
			8  & 40 & 40 & 10 & -2 & 32 & 0 & 0 & 6 & 6 & 24 & 20 & 2 & 2 & 0 & 14 & -2 & 14 & -18 & 12 & 2 & -12 & 2 & 12 & 2 & 6 & 6 & 1 & 1 & -6  \\
			9 & 56 & 38 & -7 & 5 & -3 & 1 & -7 & 0 & 2 & -7 & -7 & 2 & 2 & -4 & -21 & 7 & 4 & 3 & -2 & 0 & -5 & 2 & 10 & -4 & 0 & 0 & 0 & 0 & -2 \\
			\midrule
			n ~\backslash ~[g] &24B & 24C & 24D & 24E & 24F & 24G & 24H & 24I & 24J & 25A & 26A & 26B & 27A & 27B & 28A & 28B & 28C & 28D & 29A & 30A & 30B & 30C & 30D & 30E & 30F & 30G & 31A & 31B & 32A \\
			\midrule
			-15 &1 & 1 & 1 & 1 & 1 & 1 & 1 & 1 & 1 & 1 & 1 & 1 & 1 & 1 & 1 & 1 & 1 & 1 & 1 & 1 & 1 & 1 & 1 & 1 & 1 & 1 & 1 & 1& 1   \\
			-7 & 1 & 1 & 1 & 1 & 1 & 1 & 1 & 1 & 1 & 1 & 1 & 1 & 1 & 1 & 1 & 1 & 1 & 1 & 1 & 1 & 1 & 1 & 1 & 1 & 1 & 1 & 1 & 1 & 1  \\
			-3	& 1 & 1 & 1 & 1 & 1 & 1 & 1 & 1 & 1 & 1 & 1 & 1 & 1 & 1 & 1 & 1 & 1 & 1 & 1 & 1 & 1 & 1 & 1 & 1 & 1 & 1 & 1 & 1  & 1 	\\
			0     &-2 & -2 & -2 & -2 & -2 & -2 & -2 & -2 & -2 & -2 & -2 & -2 & -2 & -2 & -2 & -2 & -2 & -2 & -2 & -2 & -2 & -2 & -2 & -2 & -2 & -2 & -2 & -2 & -2     \\
			1 & -3 & 0 & 1 & 0 & 0 & 0 & -2 & -4 & 0 & -4 & -4 & -3 & -3 & -3 & -3 & -3 & 1 & -2 & -3 & -3 & -4 & 0 & -3 & 0 & -3 & -1 & -3 & -3 & -2 \\
			4  & -6 & 0 & 2 & 0 & -6 & 2 & -10 & -4 & 0 & -4 & -8 & -4 & -6 & -6 & -2 & -6 & -2 & -6 & -4 & 0 & -2 & -6 & 0 & 0 & -8 & -6 & -3 & -3 & -4 \\
			5  &-5 & -5 & 3 & -3 & 7 & -1 & 3 & 3 & 1 & 4 & 3 & 0 & 2 & 2 & -1 & 7 & 7 & 1 & 1 & 4 & 6 & 2 & -1 & -2 & 4 & 2 & 3 & 3 & 1 \\
			8  & 10 & 4 & -6 & 6 & -14 & 2 & 6 & 0 & -2 & 0 & 8 & 2 & 0 & 0 & 4 & -8 & -12 & 4 & 2 & -10 & -4 & 0 & 0 & 2 & 2 & -2 & -1 & -1 & 0  \\
			9 & 6 & 9 & 2 & 1 & 1 & 1 & -1 & -3 & 1 & -4 & -5 & -2 & -1 & -1 & -11 & 5 & 1 & -1 & -3 & -1 & -4 & 0 & -6 & 5 & -1 & 1 & -2 & -2 & -1 \\
			\midrule
			n ~\backslash ~[g] &32B & 33A & 33B & 34A & 35A & 35B & 36A & 36B & 36C & 36D & 38A & 39A & 39B & 39C & 39D & 40A & 40B & 40C & 40D & 41A & 42A & 42B & 42C & 42D & 44A & 44B & 45A & 46A & 46B \\
			\midrule
			-15 &1 & 1 & 1 & 1 & 1 & 1 & 1 & 1 & 1 & 1 & 1 & 1 & 1 & 1 & 1 & 1 & 1 & 1 & 1 & 1 & 1 & 1 & 1 & 1 & 1 & 1 & 1 & 1& 1   \\
			-7 & 1 & 1 & 1 & 1 & 1 & 1 & 1 & 1 & 1 & 1 & 1 & 1 & 1 & 1 & 1 & 1 & 1 & 1 & 1 & 1 & 1 & 1 & 1 & 1 & 1 & 1 & 1 & 1 & 1  \\
			-3	& 1 & 1 & 1 & 1 & 1 & 1 & 1 & 1 & 1 & 1 & 1 & 1 & 1 & 1 & 1 & 1 & 1 & 1 & 1 & 1 & 1 & 1 & 1 & 1 & 1 & 1 & 1 & 1  & 1 	\\
			0     &-2 & -2 & -2 & -2 & -2 & -2 & -2 & -2 & -2 & -2 & -2 & -2 & -2 & -2 & -2 & -2 & -2 & -2 & -2 & -2 & -2 & -2 & -2 & -2 & -2 & -2 & -2 & -2 & -2     \\
			1 & -2 & 1 & -2 & -3 & -1 & -2 & -3 & 0 & -1 & -2 & -2 & -3 & 0 & -2 & -2 & 0 & -2 & -1 & -1 & -2 & -2 & -1 & 0 & -1 & -1 & -1 & -2 & 0 & 0 \\
			4  & -4 & 2 & -4 & -4 & -6 & -4 & 0 & 0 & -6 & -2 & -6 & 0 & 0 & -3 & -3 & -2 & -2 & -4 & -4 & -2 & -4 & -2 & 0 & -2 & -4 & -4 & -2 & -3 & -3 \\
			5  &1 & 0 & -3 & 2 & 2 & 1 & 0 & 0 & 4 & 0 & 2 & 4 & -2 & 2 & 2 & 3 & 1 & 1 & 1 & 0 & -2 & -4 & 1 & -1 & 3 & 3 & 2 & 2 & 2 \\
			8  & 0 & 4 & 4 & 2 & -2 & 0 & 0 & 0 & -6 & 0 & 0 & -4 & 2 & 0 & 0 & -10 & -2 & 0 & 0 & 0 & 0 & 2 & 0 & 2 & -2 & -2 & 4 & -3 & -3  \\
			9 & -1 & 6 & 0 & -4 & -3 & -1 & -7 & 2 & 1 & -2 & 0 & -1 & 2 & -1 & -1 & 1 & -1 & 0 & 0 & -2 & -1 & 0 & 2 & 0 & 0 & 0 & -1 & 0 & 0 \\
			\bottomrule
		\end{tabular}
	\end{small}
\end{sidewaystable}

\begin{sidewaystable}
	\begin{small}
		\centering
		\caption{Coefficients of $q^n$ in the McKay--Thompson series $\cH_{[g]}(\tau)$ for the Monster, part three.}\smallskip
		\begin{tabular}{C@{ }@{\;}R@{ }R@{ }R@{ }R@{ }R@{ }R@{ }R@{ }R@{ }R@{ }R@{ }R@{ }R@{ }R@{ }R@{ }R@{ }R@{ }R@{ }R@{ }R@{ }R@{ }R@{ }R@{ }R@{ }R@{ }R@{ }R@{ }R@{ }R@{ }R@{ }R@{ }R@{ }R@{ }R@{ }R@{ }R@{ }R@{ }R}\toprule
			n ~\backslash ~[g] &46C & 46D & 47A & 47B & 48A & 50A & 51A & 52A & 52B & 54A & 55A & 56A & 56B & 56C & 57A & 59A & 59B & 60A & 60B & 60C & 60D & 60E & 60F & 62A & 62B & 66A & 66B & 68A & 69A \\
			\midrule
			-15 &1 & 1 & 1 & 1 & 1 & 1 & 1 & 1 & 1 & 1 & 1 & 1 & 1 & 1 & 1 & 1 & 1 & 1 & 1 & 1 & 1 & 1 & 1 & 1 & 1 & 1 & 1 & 1& 1   \\
			-7 & 1 & 1 & 1 & 1 & 1 & 1 & 1 & 1 & 1 & 1 & 1 & 1 & 1 & 1 & 1 & 1 & 1 & 1 & 1 & 1 & 1 & 1 & 1 & 1 & 1 & 1 & 1 & 1 & 1  \\
			-3	& 1 & 1 & 1 & 1 & 1 & 1 & 1 & 1 & 1 & 1 & 1 & 1 & 1 & 1 & 1 & 1 & 1 & 1 & 1 & 1 & 1 & 1 & 1 & 1 & 1 & 1 & 1 & 1  & 1 	\\
			0     &-2 & -2 & -2 & -2 & -2 & -2 & -2 & -2 & -2 & -2 & -2 & -2 & -2 & -2 & -2 & -2 & -2 & -2 & -2 & -2 & -2 & -2 & -2 & -2 & -2 & -2 & -2 & -2 & -2     \\
			1 & -2 & -2 & -1 & -1 & -1 & -2 & -1 & 0 & -1 & -1 & -2 & -1 & 0 & 0 & 0 & -1 & -1 & -2 & 0 & -1 & 1 & -1 & 0 & -1 & -1 & -2 & -1 & -1 & -1 \\
			4  &-3 & -3 & -3 & -3 & -2 & -2 & -2 & -4 & -2 & -4 & 0 & 0 & -2 & -2 & 0 & -2 & -2 & 0 & -2 & -2 & 0 & -2 & 0 & -3 & -3 & 0 & 0 & -2 & -1 \\
			5  & 2 & 2 & 1 & 1 & -1 & 0 & -1 & 3 & 0 & 2 & 0 & -1 & 3 & 3 & 0 & 1 & 1 & 0 & 2 & 2 & 4 & 1 & 0 & 1 & 1 & 1 & 0 & 2 & -1\\
			8  & -1 & -1 & -1 & -1 & 2 & 0 & 0 & -4 & 0 & 0 & -4 & -2 & -4 & -4 & 0 & 0 & 0 & -2 & -4 & -4 & -4 & 0 & 0 & -3 & -3 & 0 & 0 & -4 & 1  \\
			9 & 0 & 0 & -1 & -1 & 0 & -2 & 0 & -1 & 0 & -1 & -2 & -1 & 1 & 1 & 2 & 0 & 0 & -2 & 0 & 1 & 3 & 0 & 1 & 0 & 0 & 0 & 0 & 2 & 0\\
			\midrule
			n ~\backslash ~[g] &69B & 70A & 70B & 71A & 71B & 78A & 78B & 78C & 84A & 84B & 84C & 87A & 87B & 88A & 88B & 92A & 92B & 93A & 93B & 94A & 94B & 95A & 95B & 104A & 104B & 105A & 110A & 119A & 119B\\
			\midrule
			-15 &1 & 1 & 1 & 1 & 1 & 1 & 1 & 1 & 1 & 1 & 1 & 1 & 1 & 1 & 1 & 1 & 1 & 1 & 1 & 1 & 1 & 1 & 1 & 1 & 1 & 1 & 1 & 1& 1   \\
			-7 & 1 & 1 & 1 & 1 & 1 & 1 & 1 & 1 & 1 & 1 & 1 & 1 & 1 & 1 & 1 & 1 & 1 & 1 & 1 & 1 & 1 & 1 & 1 & 1 & 1 & 1 & 1 & 1 & 1  \\
			-3	& 1 & 1 & 1 & 1 & 1 & 1 & 1 & 1 & 1 & 1 & 1 & 1 & 1 & 1 & 1 & 1 & 1 & 1 & 1 & 1 & 1 & 1 & 1 & 1 & 1 & 1 & 1 & 1  & 1 	\\
			0     &-2 & -2 & -2 & -2 & -2 & -2 & -2 & -2 & -2 & -2 & -2 & -2 & -2 & -2 & -2 & -2 & -2 & -2 & -2 & -2 & -2 & -2 & -2 & -2 & -2 & -2 & -2 & -2 & -2     \\
			1 &-1 & -1 & 0 & -1 & -1 & -1 & 0 & 0 & 0 & 1 & 0 & 0 & 0 & -1 & -1 & 0 & 0 & 0 & 0 & -1 & -1 & -1 & -1 & 0 & 0 & -1 & 0 & 0 & 0\\
			4  & -1 & -2 & -2 & -1 & -1 & -2 & -1 & -1 & -2 & 0 & 0 & -1 & -1 & 0 & 0 & -1 & -1 & 0 & 0 & -1 & -1 & -1 & -1 & 0 & 0 & 0 & -2 & -1 & -1 \\
			5  &-1 & 2 & 3 & 0 & 0 & 0 & 0 & 0 & 2 & 4 & 1 & 1 & 1 & -1 & -1 & 2 & 2 & 0 & 0 & 1 & 1 & 1 & 1 & 1 & 1 & -1 & 2 & 2 & 2 \\
			8  &1 & 2 & 0 & -2 & -2 & -4 & -4 & -4 & -2 & -2 & -2 & -1 & -1 & -2 & -2 & -3 & -3 & -1 & -1 & -1 & -1 & 0 & 0 & -4 & -4 & -2 & -4 & -1 & -1  \\
			9 & 0 & 1 & 1 & 0 & 0 & 1 & 1 & 1 & 1 & 2 & 2 & 0 & 0 & 0 & 0 & 0 & 0 & 1 & 1 & 1 & 1 & 1 & 1 & 1 & 1 & 0 & 0 & 0 & 0\\
			\bottomrule
		\end{tabular}
	\end{small}
\end{sidewaystable}

\begin{sidewaystable}
	\begin{small}
		\centering
		\caption{Coefficients of $q^n$ in the McKay--Thompson series $\cG_{[g]}(\tau)$ for the Baby monster, part one.}\label{mtM}\smallskip
	\resizebox{!}{56pt}{	\begin{tabular}{C@{ }@{\;}R@{ }R@{ }R@{ }R@{ }R@{ }R@{ }R@{ }R@{ }R@{ }R@{ }R@{ }R@{ }R@{ }R@{ }R@{ }R@{ }R@{ }R@{ }R@{ }R@{ }R@{ }R@{ }R@{ }R@{ }R@{ }R@{ }R@{ }R@{ }R@{ }R@{ }R}\toprule
			n ~\backslash ~[g] &
	1A & 2A & 2B & 2C & 2D & 3A & 3B & 4A & 4B & 4C & 4D & 4E & 4F & 4G & 4H & 4I & 4J & 5A & 5B & 6A & 6B & 6C & 6D & 6E   \\
			\midrule
			-7 & 1 & 1 & 1 & 1 & 1 & 1 & 1 & 1 & 1 & 1 & 1 & 1 & 1 & 1 & 1 & 1 & 1 & 1 & 1 & 1 & 1 & 1 & 1 & 1  \\
			-3	& 1 & 1 & 1 & 1 & 1 & 1 & 1 & 1 & 1 & 1 & 1 & 1 & 1 & 1 & 1 & 1 & 1 & 1 & 1 & 1 & 1 & 1 & 1 & 1  	\\
			0     &-2 & -2 & -2 & -2 & -2 & -2 & -2 & -2 & -2 & -2 & -2 & -2 & -2 & -2 & -2 & -2 & -2 & -2 & -2 & -2 & -2 & -2 & -2 & -2 \\
			1 &-4371 & 493 & -275 & 53 & -19 & -78 & 3 & 77 & -51 & -19 & 21 & -35 & 13 & -11 & -3 & 1 & 5 & -21 & 4 & 34 & -20 & -14 & 7 & -4  \\
			4  & 8315004 & 141436 & 24700 & -1668 & -1924 & 1650 & 30 & 3452 & 380 & 124 & 124 & 252 & -132 & -132 & -4 & 52 & -4 & 54 & 4 & 658 & 10 & -14 & 10 & -22 \\
			5  & -52842475 & 565269 & -85995 & 9237 & 21 & -4069 & -181 & 9749 & -491 & 21 & 533 & -491 & 533 & 21 & 21 & -51 & 21 & -75 & 25 & 1563 & -3 & 27 & 51 & -3  \\
			8  & 5736480000 & 18517248 & 1752320 & 46848 & 44288 & 34890 & 870 & 102144 & 5888 & 3328 & 3328 & 2816 & 768 & 768 & 256 & -56 & 256 & 350 & 0 & 10314 & 378 & 74 & 162 & 122  \\
		\end{tabular}}
		\resizebox{!}{106pt}{
		\begin{tabular}{C@{ }@{\;}R@{ }R@{ }R@{ }R@{ }R@{ }R@{ }R@{ }R@{ }R@{ }R@{ }R@{ }R@{ }R@{ }R@{ }R@{ }R@{ }R@{ }R@{ }R@{ }R@{ }R@{ }R@{ }R@{ }R@{ }R@{ }R@{ }R@{ }R@{ }R@{ }R@{ }R@{ }R@{ }R@{ }R@{ }R@{ }R@{ }R@{ }R@{ }R@{ }R@{ }R@{ }R}\midrule
			n ~\backslash ~[g] &6F & 6G & 6H & 6I & 6J & 6K & 7A & 8A & 8B & 8C & 8D & 8E & 8F & 8G & 8H & 8I & 8J & 8K & 8L & 8M & 8N & 9A & 9B & 10A & 10B & 10C & 10D & 10E & 10F & 11A & 12A & 12B\\
			\midrule
			-7  & 1 & 1 & 1 & 1 & 1 & 1 & 1 & 1 & 1 & 1 & 1 & 1 & 1 & 1 & 1 & 1 & 1 & 1 & 1 & 1 & 1 & 1 & 1 & 1 & 1 & 1 & 1 & 1 & 1 & 1 & 1 & 1 \\
			-3	&1 & 1 & 1 & 1 & 1 & 1 & 1 & 1 & 1 & 1 & 1 & 1 & 1 & 1 & 1 & 1 & 1 & 1 & 1 & 1 & 1 & 1 & 1 & 1 & 1 & 1 & 1 & 1 & 1 & 1 & 1 & 1	\\
			0     &-2 & -2 & -2 & -2 & -2 & -2 & -2 & -2 & -2 & -2 & -2 & -2 & -2 & -2 & -2 & -2 & -2 & -2 & -2 & -2 & -2 & -2 & -2 & -2 & -2 & -2 & -2 & -2 & -2 & -2 & -2 & -2   \\
			1 &8 & -5 & 2 & -13 & -1 & -1 & -10 & 21 & -7 & 9 & -11 & 1 & 5 & 5 & -3 & 1 & -3 & -7 & 1 & 1 & 1 & 3 & -6 & -7 & -5 & 3 & 0 & 1 & -4 & -4 & 12 & 3 \\
			4  &-30 & 22 & -46 & 14 & 2 & 6 & -2 & 252 & 12 & 12 & -4 & -52 & -4 & -20 & -4 & 12 & -20 & -4 & -12 & -4 & 4 & 3 & -6 & -14 & -10 & -18 & 0 & -14 & -4 & -6 & 74 & 2  \\
			5  &165 & 27 & 27 & -21 & 3 & 3 & 7 & 533 & -11 & 117 & 21 & 53 & 21 & 69 & 21 & -11 & 5 & 5 & 13 & 5 & -3 & 8 & 8 & 19 & 5 & 37 & 5 & 11 & 1 & 7 & 157 & -5  \\
			8  &138 & -34 & 74 & 86 & 26 & 30 & 20 & 2560 & 192 & 192 & 0 & 64 & 0 & 16 & 0 & 0 & 16 & 16 & 8 & 16 & -8 & 6 & 6 & -2 & 30 & -2 & 20 & -2 & 8 & -2 & 506 & 2  \\
			\midrule
			n ~\backslash ~[g]&12C & 12D & 12E & 12F & 12G & 12H & 12I & 12J & 12K & 12L & 12M & 12N & 12O & 12P & 12Q & 12R & 12S & 12T & 13A & 14A & 14B & 14C & 14D & 14E & 15A & 15B & 16A & 16B & 16C & 16D & 16E & 16F\\
			\midrule
			-7  & 1 & 1 & 1 & 1 & 1 & 1 & 1 & 1 & 1 & 1 & 1 & 1 & 1 & 1 & 1 & 1 & 1 & 1 & 1 & 1 & 1 & 1 & 1 & 1 & 1 & 1 & 1 & 1 & 1 & 1 & 1 & 1 \\
			-3	&1 & 1 & 1 & 1 & 1 & 1 & 1 & 1 & 1 & 1 & 1 & 1 & 1 & 1 & 1 & 1 & 1 & 1 & 1 & 1 & 1 & 1 & 1 & 1 & 1 & 1 & 1 & 1 & 1 & 1 & 1 & 1	\\
			0     &-2 & -2 & -2 & -2 & -2 & -2 & -2 & -2 & -2 & -2 & -2 & -2 & -2 & -2 & -2 & -2 & -2 & -2 & -2 & -2 & -2 & -2 & -2 & -2 & -2 & -2 & -2 & -2 & -2 & -2 & -2 & -2   \\
			1 &-4 & -6 & 6 & -5 & 2 & 0 & 4 & -2 & 3 & -2 & -1 & 1 & -5 & 0 & 3 & -1 & -3 & 1 & -3 & 10 & -4 & 4 & -2 & 2 & -3 & -2 & -3 & 5 & 1 & -3 & 1 & -3\\
			4  &-22 & 2 & 10 & 6 & -14 & -14 & -6 & -6 & -2 & -6 & -2 & 0 & -6 & 2 & -10 & -10 & -4 & -2 & -4 & 50 & -6 & -2 & -10 & -6 & 0 & -5 & 4 & 4 & -4 & -4 & -12 & -4 \\
			5  &29 & -5 & 35 & -5 & 27 & 5 & 29 & -5 & 11 & 3 & 3 & -5 & 11 & -3 & 27 & 3 & 3 & 3 & 3 & 103 & 5 & 25 & 7 & 7 & 6 & 4 & -3 & 29 & 5 & 5 & 13 & 5 \\
			8  &-6 & 74 & 58 & 62 & 10 & 10 & -6 & 26 & -2 & -6 & -2 & 8 & -6 & 10 & -2 & -2 & 4 & -2 & 4 & 288 & 8 & 4 & -4 & -8 & -10 & 0 & 24 & 24 & 0 & 0 & -8 & 0 \\
			\bottomrule
		\end{tabular}}
	\end{small}
\end{sidewaystable}

\begin{sidewaystable}
	\begin{small}
		\centering
		\caption{Coefficients of $q^n$ in the McKay--Thompson series $\cG_{[g]}(\tau)$ for the Baby monster, part two.}\smallskip
		\resizebox{!}{157pt}{\begin{tabular}{C@{ }@{\;}R@{ }R@{ }R@{ }R@{ }R@{ }R@{ }R@{ }R@{ }R@{ }R@{ }R@{ }R@{ }R@{ }R@{ }R@{ }R@{ }R@{ }R@{ }R@{ }R@{ }R@{ }R@{ }R@{ }R@{ }R@{ }R@{ }R@{ }R@{ }R@{ }R@{ }R@{ }R@{ }R@{ }R@{ }R@{ }R@{ }R@{ }R@{ }R@{ }R@{ }R@{ }R}\midrule
			n ~\backslash ~[g]&16G & 16H & 17A & 18A & 18B & 18C & 18D & 18E & 18F & 19A & 20A & 20B & 20C & 20D & 20E & 20F & 20G & 20H & 20I & 20J & 21A & 22A & 22B & 23A & 23B & 24A & 24B & 24C & 24D & 24E & 24F & 24G\\
			\midrule
			-7  & 1 & 1 & 1 & 1 & 1 & 1 & 1 & 1 & 1 & 1 & 1 & 1 & 1 & 1 & 1 & 1 & 1 & 1 & 1 & 1 & 1 & 1 & 1 & 1 & 1 & 1 & 1 & 1 & 1 & 1 & 1 & 1 \\
			-3	&1 & 1 & 1 & 1 & 1 & 1 & 1 & 1 & 1 & 1 & 1 & 1 & 1 & 1 & 1 & 1 & 1 & 1 & 1 & 1 & 1 & 1 & 1 & 1 & 1 & 1 & 1 & 1 & 1 & 1 & 1 & 1	\\
			0     &-2 & -2 & -2 & -2 & -2 & -2 & -2 & -2 & -2 & -2 & -2 & -2 & -2 & -2 & -2 & -2 & -2 & -2 & -2 & -2 & -2 & -2 & -2 & -2 & -2 & -2 & -2 & -2 & -2 & -2 & -2 & -2   \\
			1 &-1 & -1 & -2 & -2 & -2 & -2 & 2 & -1 & 2 & -1 & 7 & 2 & -3 & -1 & 1 & -1 & 0 & 1 & -2 & 0 & -1 & -2 & 0 & -1 & -1 & 2 & 0 & -2 & 4 & 0 & 3 & -2 \\
			4  &-4 & -4 & -2 & -8 & -8 & -2 & -4 & -1 & -6 & -4 & 22 & 2 & 2 & -10 & -6 & -2 & -8 & 2 & -2 & -4 & -2 & -2 & -6 & -2 & -2 & -6 & -6 & 2 & 2 & -6 & 0 & 2  \\
			5  &1 & 1 & 2 & 6 & 6 & 0 & 6 & 0 & 12 & 2 & 49 & -1 & -1 & 9 & 3 & 1 & 9 & -1 & 3 & 1 & -2 & 1 & 3 & 2 & 2 & 7 & 9 & -1 & 17 & 5 & 11 & 3  \\
			8  &-4 & 4 & 0 & 0 & 0 & 2 & -4 & 2 & -6 & 0 & 94 & -6 & 14 & -2 & -2 & -2 & -4 & -6 & -2 & -4 & 2 & 2 & -2 & -2 & -2 & -6 & -6 & 10 & 10 & 10 & 4 & -6 \\
			\midrule
			n ~\backslash ~[g]&24H & 24I & 24J & 24K & 24L & 24M & 24N & 25A & 26A & 26B & 27A & 28A & 28B & 28C & 28D & 28E & 30A & 30B & 30C & 30D & 30E & 30F & 30G & 30H & 31A & 31B & 32A & 32B & 32C & 32D & 33A & 34A\\
			\midrule
			-7  & 1 & 1 & 1 & 1 & 1 & 1 & 1 & 1 & 1 & 1 & 1 & 1 & 1 & 1 & 1 & 1 & 1 & 1 & 1 & 1 & 1 & 1 & 1 & 1 & 1 & 1 & 1 & 1 & 1 & 1 & 1 & 1 \\
			-3	&1 & 1 & 1 & 1 & 1 & 1 & 1 & 1 & 1 & 1 & 1 & 1 & 1 & 1 & 1 & 1 & 1 & 1 & 1 & 1 & 1 & 1 & 1 & 1 & 1 & 1 & 1 & 1 & 1 & 1 & 1 & 1	\\
			0     &-2 & -2 & -2 & -2 & -2 & -2 & -2 & -2 & -2 & -2 & -2 & -2 & -2 & -2 & -2 & -2 & -2 & -2 & -2 & -2 & -2 & -2 & -2 & -2 & -2 & -2 & -2 & -2 & -2 & -2 & -2 & -2   \\
			1 &-1 & 2 & 0 & -1 & -1 & 1 & 1 & -1 & 1 & -1 & 0 & -2 & 2 & 0 & 0 & 0 & 5 & -1 & 3 & 1 & 1 & 2 & 0 & 0 & 0 & 0 & -1 & -1 & 1 & 1 & -1 & 2\\
			4  & -2 & -2 & -2 & -4 & 2 & -6 & -2 & -1 & -4 & -4 & -3 & 2 & -2 & -6 & -2 & 0 & 10 & -2 & 0 & -4 & -2 & -1 & -3 & -3 & -2 & -2 & 0 & 0 & -4 & -4 & 0 & -2\\
			5  &3 & 9 & -1 & 3 & -1 & 7 & 3 & 0 & 7 & 3 & 2 & -1 & 7 & 5 & 1 & -1 & 22 & -2 & 10 & 2 & 2 & 4 & 2 & 2 & 1 & 1 & 1 & 1 & 5 & 5 & 1 & 6 \\
			8  &-2 & -2 & -2 & 0 & -2 & -6 & -2 & 0 & -4 & -4 & 0 & 8 & -4 & 0 & -4 & 2 & 28 & 4 & -2 & -6 & -8 & -4 & -4 & -4 & -2 & -2 & 0 & 0 & -4 & -4 & -2 & -4 \\
			\midrule
			n ~\backslash ~[g] &34B & 34C & 35A & 36A & 36B & 36C & 38A & 39A & 40A & 40B & 40C & 40D & 40E & 42A & 42B & 42C & 44A & 46A & 46B & 47A & 47B & 48A & 48B & 52A & 55A & 56A & 56B & 60A & 60B & 60C & 66A & 70A\\
			\midrule
			-7  & 1 & 1 & 1 & 1 & 1 & 1 & 1 & 1 & 1 & 1 & 1 & 1 & 1 & 1 & 1 & 1 & 1 & 1 & 1 & 1 & 1 & 1 & 1 & 1 & 1 & 1 & 1 & 1 & 1 & 1 & 1 & 1 \\
			-3	&1 & 1 & 1 & 1 & 1 & 1 & 1 & 1 & 1 & 1 & 1 & 1 & 1 & 1 & 1 & 1 & 1 & 1 & 1 & 1 & 1 & 1 & 1 & 1 & 1 & 1 & 1 & 1 & 1 & 1 & 1 & 1	\\
			0     &-2 & -2 & -2 & -2 & -2 & -2 & -2 & -2 & -2 & -2 & -2 & -2 & -2 & -2 & -2 & -2 & -2 & -2 & -2 & -2 & -2 & -2 & -2 & -2 & -2 & -2 & -2 & -2 & -2 & -2 & -2 & -2   \\
			1 &0 & 0 & 0 & 0 & 0 & 1 & -1 & 0 & 3 & -1 & 1 & 1 & 0 & 1 & -1 & 1 & 0 & 1 & 1 & 0 & 0 & 0 & 2 & 1 & 1 & 0 & 0 & 1 & 1 & 0 & 1 & 0\\
			4  & -4 & -4 & -2 & -4 & -2 & -3 & 0 & -1 & 2 & 2 & -2 & -2 & 0 & -4 & 0 & -2 & -2 & -2 & -2 & -1 & -1 & -2 & -2 & 0 & -1 & 0 & 0 & -2 & 0 & 1 & -2 & 0\\
			5  &2 & 2 & 2 & 4 & 2 & 4 & 0 & 0 & 9 & -3 & 3 & 3 & -1 & 4 & 2 & 4 & 3 & 2 & 2 & 1 & 1 & 3 & 5 & 1 & 2 & 1 & 1 & 4 & 0 & 0 & 1 & -2 \\
			8  &-2 & -2 & 0 & -4 & -2 & -4 & 0 & -2 & 2 & 2 & -6 & -2 & -4 & 0 & -4 & -2 & -2 & -4 & -4 & -2 & -2 & -6 & -6 & -4 & -2 & -2 & -2 & 4 & -2 & 2 & -4 & -2\\
			\bottomrule
		\end{tabular}}
	\end{small}
\end{sidewaystable}

\end{document}